\documentstyle [11pt]{article}

\newfam\msbfam
\font\twelmsb=msbm10 at 12pt
\font\tenmsb=msbm10
\font\sevenmsb=msbm10 at 7pt
\font\fivemsb=msbm10 at 5pt

\textfont\msbfam=\twelmsb
\scriptfont\msbfam=\sevenmsb
\scriptscriptfont\msbfam=\fivemsb

\def\Bbb{\fam\msbfam\tenmsb}

\def\C{{\Bbb C}}

\def\P{{\Bbb P}}

\def\R{{\Bbb R}}
\def\Z{{\Bbb Z}}

\def\KK{{\cal K}}

\def\OO{{\cal O}}

\def\WW{{\cal W}}

\def\map{\longrightarrow}
\def\textmap#1{\mathop{\vbox{\ialign{
                                ##\crcr
    ${\scriptstyle\hfil\;\;#1\;\;\hfil}$\crcr
    \noalign{\kern-1pt\nointerlineskip}
    \rightarrowfill\crcr}}\;}}

\def\textlmap#1{\mathop{\vbox{\ialign{
                                ##\crcr
    ${\scriptstyle\hfil\;\;#1\;\;\hfil}$\crcr
    \noalign{\kern-1pt\nointerlineskip}
    \leftarrowfill\crcr}}\;}}

\newfam\meuffam
\font\tenmeuf=eufm10
\font\sevenmeuf=eufm7
\font\fivemeuf=eufm5

\textfont\meuffam=\tenmeuf
\scriptfont\meuffam=\sevenmeuf
\scriptscriptfont\meuffam=\fivemeuf

\def\germ{\fam\meuffam\tenmeuf}

\def\picture#1by#2(#3){
\vbox to #2 {
  \hrule width #1 height 0pt depth 0pt \vfill \special{picture #3}}
}

\def\scaledpicture#1by#2(#3scaled#4){{
\dimen0=#1  \dimen1=#2
\divide\dimen0 by 1000 \multiply\dimen0 by #4
\divide\dimen1 by 1000 \multiply\dimen1 by #4
\picture \dimen0 by \dimen1 (#3 scaled #4)}}
\def\dfigure#1by#2(#3scaled#4offset#5:#6)
  {\medskip
   \vglue 2mm minus 2mm
   $$
     \hbox{
       \hglue#5
       {\scaledpicture #1 by #2 (#3 scaled #4)}
     }
   $$
   \par\goodbreak
   \vglue 2mm minus 2mm
   \medskip}

\begin{document}
\def\Pr{{\rm Pr}}
\def\tr{{\rm Tr}}
\def\End{{\rm End}}
\def\Pic{{\rm Pic}}
\def\NS{{\rm NS}}
\def\deg{{\rm deg}}
\def\Hom{{\rm Hom}}
\def\Herm{{\rm Herm}}
\def\Vol{{\rm Vol}}
\def\pf{{\bf Proof: }}
\def\id{{\rm id}}
\def\i{{\germ i}}
\def\im{{\rm im}}
\def\rk{{\rm rk}}
\def\Spin{{\rm Spin}}
\def\h{{\bf H}}
\def\dv{\bar\partial}
\def\dva{\bar\partial_A}
\def\da{\partial_A}
\def\p{\partial\bar\partial}
\def\pa{\partial_A\bar\partial_A}
\def\Dr{\hskip 4pt{\not}{D}}
\newtheorem{sz}{Satz}
\newtheorem{szfr}{Satzfr}
\newtheorem{th}[sz]{Theorem}
\newtheorem{thfr}[szfr]{Th\'eor\`eme}
\newtheorem{pr}[sz]{Proposition}
\newtheorem{re}[sz]{Remark}
\newtheorem{co}[sz]{Corollary}
\newtheorem{cofr}[szfr]{Corollaire}
\newtheorem{dt}[sz]{Definition}
\newtheorem{lm}[sz]{Lemma}

\noindent G\'eom\'etrie diff\'erentielle/{\sl Differential Geometry}
\\ \\ \\
\centerline{\Large{\bf Les Invariants de Seiberg-Witten et}} \\
\vspace{5mm}
\centerline{\Large{\bf la Conjecture de Van De Ven}}
\vspace{5mm}
\centerline{\large Christian Okonek\footnotemark  \ \ et \ \ Andrei
Teleman$^1$}\\ \footnotetext[1]{Partially supported by: AGE-Algebraic
Geometry in
Europe, contract No ERBCHRXCT940557 (BBW 93.0187), and by  SNF, nr.
21-36111.92}
{\footnotesize{\sl R\'esum\'e --}} {\footnotesize Nous pr\'esentons une
d\'emonstration tr\`es courte du fait que la rationalit\'e des surfaces
complexes est une propri\'et\'e ne d\'ependant que de la structure
diff\'erentielle. La preuve utilise les nouveaux invariants introduits
par Seiberg-Witten.}
\vspace{3mm}\\
\centerline{{\bf Seiberg-Witten Invariants and the Van De Ven Conjecture}}
{\footnotesize{\sl Abstract --}} {\footnotesize We give a short proof for
the fact that rationality of complex surfaces is a property depending only
on the differential structure. Our proof uses the new Seiberg-Witten
invariants.}
\vspace{3mm}\\
{\sl Abridged English Version} -- We prove the following theorem [8]
\begin{th}
A complex surface which is diffeomeorphic to a rational surface is rational.
\end{th}
This result has been announced by R.\ Friedman and Z.\ Qin [4].
Whereas their proof uses Donaldson theory and vector bundles techniques, our
proof uses the new Seiberg-Witten invariants [9], and the interpretation of
these invariants in terms of stable pairs [7].

Combining the theorem above with the results of [3], one obtains a proof
of the Van de Ven conjecture [8], [10]:
\begin{co}
The Kodaira dimension of a complex surface is a differential invariant.
\end{co}

{\bf Proof:} (of the Theorem) It suffices to prove the theorem for
algebraic surfaces [1].
Let $X$ be an algebraic surface of non-negative Kodaira dimension, with
$\pi_1(X)=\{1\}$ and $p_g(X)=0$. We may suppose that $X$ is the blow up
in $k$ {\sl distinct} points of its minimal model $X_{\min}$. Denote the
contraction to the minimal model by $\sigma:X\map X_{\min}$, and the
exceptional divisor by $E=\sum\limits_{i=1}^kE_i$.

Fix an ample divisor $H_{\min}$ on $X_{\min}$, a sufficiently large
integer $n$, and let $H_n:=\sigma^*(n H_{\min})-E$ be the associated
polarization of $X$.

For every subset $I\subset\{1,\dots,k\}$ we put
$E_I:=\sum\limits_{i\in I} E_i$, and $c_I:=2[E_I]-[K_X]$, where $K_X$ is
a canonical divisor. Clearly $c_I=[E_I]-[E_{\bar I}]-\sigma^*([K_{\min}])$,
where $\bar I$ denotes the complement of $I$ in $\{1,\dots,k\}$. The
cohomology classes $c_I$ are almost canonical classes in the sense of [7].
Now choose a K\"ahler metric $g_n$ on $X$ with K\"ahler class
$[\omega_{g_n}]=c_1(\OO_X(H_n))$. Since $[\omega_{g_n}]\cdot c_I<0$ for
sufficiently large $n$, the main result of [7] identifies the Seiberg-Witten
moduli space $\WW_X^{g_n}(c_I)$ with the union of all complete linear
systems $|D|$ corresponding to effective divisors $D$ on $X$ with
$c_1(\OO_X(2D-K_X))=c_I$.

Since $H^2(X,\Z)$ has no 2-torsion, and $q(X)=0$, there is only one such
divisor, $D=E_I$. Furthermore, from
$h^1(\OO_X(E_I)|_{E_I})=0$, and the smoothness criterion in [7], we obtain:
$$\WW_X^{g_n}(c_I)=\{E_I\}, $$
i.e. $\WW_X^{g_n}(c_I)$ consists of a single smooth point.
The corresponding Seiberg-Witten invariants are therefore odd:
$n_{c_I}^{g_n}=\pm 1$.

Consider now the positive cone  $\KK:=\{u\in H^2_{\rm DR}(X)|\ u^2>0\}$;
using the Hodge index theorem, the fact that $K_{\min}$ is cohomologically
nontrivial, and $K_{\min}^2\geq 0$, we see that
$\KK$ splits as a disjoint union of two components  $\KK_{\pm}:=\{u\in\KK|\
\pm u\cdot\sigma^*(K_{\min})>0\}$. Clearly $[\omega_{g_n}]$ belongs
to $\KK_+$.

Let $g$ be an {\sl arbitrary} Riemannian metric on $X$, and let $\omega_g$ be
a $g$-selfdual closed 2-form on $X$ such that $[\omega_g]\in\KK_+$.

For a fixed $I\subset\{1,\dots,k\}$, we denote by $c_I^{\bot}\subset\KK_+$ the
wall associated with $c_I$, i.e. the subset of classes $u$ with $u\cdot
c_I=0$.\\ \\
{\bf Claim:} The rays $\R_{>0}[\omega_g]$,  $\R_{>0}[\omega_{g_n}]$ belong
either to the same component of $\KK_+\setminus c_I^{\bot}$ or to the same
component of   $\KK_+\setminus c_{\bar I}^{\bot}$. \\

Indeed, since $[\omega_{g_n}]\cdot c_I<0$ and $[\omega_{g_n}]\cdot c_{\bar
I}<0$,
we just have to exclude that
$$[\omega_{g}]\cdot c_I\geq 0 \ \ \ \ {\rm and}\ \ \ \ [\omega_{g}]\cdot
c_{\bar I}\geq 0 .\eqno{(*)}$$

Write $[\omega_g]=\sum\limits_{i=1}^{k}\lambda_i[E_i]+\sigma^*[\omega]$,
for some class
$[\omega]\in H^2_{\rm DR}(X_{\min})$; then
$[\omega]^2>\sum\limits_{i=1}^k\lambda_i^2$, and $[\omega]\cdot
K_{\min}>0$, since
$\omega_g$ was chosen such that its cohomology class belongs to $\KK_+$. The
inequalities $(*)$ can now be written as
$$-\sum\limits_{i\in I}\lambda_i+\sum\limits_{j\in\bar I}\lambda_j-
[\omega]\cdot K_{\min}\geq 0 \ \ \ {\rm and}\ \ \
-\sum\limits_{j\in\bar I}\lambda_j+\sum\limits_{i\in I}\lambda_i-
[\omega]\cdot K_{\min}\geq 0,$$
and we obtain the contradiction $[\omega]\cdot K_{\min}\leq 0$. This proves
the claim.

We know already that the mod 2 Seiberg-Witten invariants $n^{g_n}_{c_I}$(mod 2)
and  $n^{g_n}_{c_{\bar I}}$(mod 2) are nontrivial for the special metric $g_n$.
Since the invariants $n^{g}_{c_{I}}$(mod 2) and $n^{g}_{c_{\bar I}}$(mod 2)
depend only on the chamber of the ray $\R_{>0}[\omega_g]$ with respect to the
wall $c_{I}^{\bot}$, respectively $c_{\bar I}^{\bot}$ (see [9], [6]), at
least one of the invariants associated with the metric $g$ must be
non-zero, too.

But any rational surface admits a Hodge metric with positive total scalar
curvature [5], and with respect to such a metric {\sl all} Seiberg-Witten
invariants are trivial [7].
\vspace{0mm}\\
\noindent\hbox to 5cm{\hrulefill}
\vspace{0mm}\\
Soit $(X,g)$ une vari\'et\'e riemannienne simplement connexe, compacte et
orien\-t\'ee,  $c\in H^2(X,\Z)$ un el\'ement caract\'eristique, et $L_c$
le $S^1$-fibr\'e associ\'e. Nous d\'esignons par
$\Sigma_c=\Sigma^+_c\oplus\Sigma^-_c$ le fibr\'e vectoriel de spineurs
d\'efini par la $\Spin^c$-structure associ\'ee \`a $c$. Si
$a\in{\cal A}(L_c)$ est une $S^1$-connexion sur $L_c$, on d\'esigne
par $\Dr_a$ l'op\'erateur de Dirac correspondant.

Les \'equations de Seiberg-Witten perturb\'ees $(SW_{\mu})$ pour une
paire $(a,\Psi)\in{\cal A}(L_c)\times A^0(\Sigma^+_c)$ s'\'ecrivent [9]:
$$\left\{\begin{array}{lll} \Dr_a\Psi&=&0\\
\Gamma(F_a^+ +i\mu)&=&2(\Psi\otimes\bar\Psi)_0\ ,
\end{array}\right.\eqno{(SW_{\mu})}$$
o\`u $\Gamma:\Lambda^2_+\otimes\C\map\End_0(\Sigma^+_c)$ est l'isomorphisme
d\'efini par la $\Spin^c$-structure associ\'ee \`a $c$, et $\mu\in A^2_+$ un
param\`etre qui doit \^etre interpr\'et\'e comme une perturbation des
\'equations $(SW_0)$. Le terme $(\Psi\otimes\bar\Psi)_0$ dans la seconde
\'equation d\'esigne la composante \`a trace nulle de l'endomorphisme
hermitien $(\Psi\otimes\bar\Psi)$. Le groupe de jauge
${\cal C}^{\infty}(X,S^1)$ est ab\'elien, et op\`ere sur l'ensemble des
solutions de $(SW_{\mu})$. Soit ${\cal W}^{g,\mu}_X(c)$ l'\'espace des
modules des solutions, muni de la structure d'\'espace analytique r\'eel
naturelle. La {\sl dimension pr\'esum\'ee} de cet espace est
$w_c:=\frac{1}{4}(c^2-2e(X)-3\sigma(X))$.

Supposons que le repr\'esentant $g$-harmonique de $c$ ne soit pas
anti-auto-dual. On dit alors que $g$ est $c$-bonne. Alors, pour toute petite
perturbation $\mu$ suffisamment g\'en\'erale, ${\cal W}^{g,\mu}_X(c)$ est
lisse, compact et a la dimension pr\'esum\'ee. La classe de cobordisme d'un
tel espace de modules ne d\'epend pas de la perturbation $\mu$ [6]. Pour une
classe $c$ telle que $w_c=0$, soit $n^g_c$ (mod 2) le nombre de points
mod 2 d'un \'espace de modules ${\cal W}^{g,\mu}_X(c)$ qui a les
propriet\'es mention\'ees. Ce nombre est aussi ind\'ependant de la
m\'etrique $g$ si on suppose que $b_+\geq 2$ ou que
$c^2\geq$ et $c\ne 0$, donc il d\'efinit dans ces cas un invariant
diff\'erentiel de la vari\'et\'e.

Dans le cas $b_+=1$, $b_2>1$, $n^g_c$ d\'epend
de $g$ de la fa\c con suivante:

Soit $\KK_+$ l'une des deux composantes connexes du c\^one  \linebreak
$\KK:=\{u\in
H^2_{\rm DR}(X)\ |\ u^2>0\}$. Pour chaque m\'etrique riemannienne
$g$ on choisit une forme $g$-harmonique auto-duale $\omega_g$ telle que
$[\omega_g]\in \KK_+$. Pour deux m\'etriques $g_0$, $g_1$ on a alors
$n_c^{g_0}=n_c^{g_1}$ si et seulement si les demi-droites
$\R_{>0}[\omega_{g_0}]$, $\R_{>0}[\omega_{g_1}]$ se trouvent du m\^eme
c\^ot\'e de l'hyperplan $c^{\bot}:=\{u\in H^2_{\rm DR} \ |\ c\cdot u=0\}$.
\dfigure 80mm by 160mm (kegel scaled 500 offset 1mm:)

Soit maintenant  $(X,J,g)$ une surface complexe dot\'ee d'une m\'etrique
k\"ahlerienne. Les fibr\'es de spineurs de la structure spinorielle
canonique sont
$\Sigma^+=\Lambda^{00}\oplus\Lambda^{02}$, $\Sigma^-=\Lambda^{01}$ [7],
et pour une classe caract\'eristique arbitraire $c$ on a
$\Sigma_c^{\pm}=\Sigma^{\pm}\otimes M$, o\`u $M$ est le fibr\'e hermitien en
droites d\'efini  par la relation
$2c_1(M)+c_1(K_X^{\vee})=c$.

\begin{thfr} {\rm [7], [9]}
Soit $(X,J,g)$ une surface k\"ahlerienne simplement con\-ne\-xe, $\omega_g$
la forme de K\"ahler associ\'ee et $c\in H^2(X,\Z)$ un el\'ement
caract\'eristique tel que $\pm c\cup[\omega_g]<0$.
\hfill{\break}
1. \ Si $c\ \not\in\ \NS(X)$\ ,  ${\cal W}_X^g(c)=\emptyset$ .\hfill{\break}
2. Supposons que $c\in\NS(X)$. Alors il existe un isomorphisme
r\'eel-analytique
naturel  ${\cal W}_X^g(c)\simeq\P(H^0(X,{\cal M}))$, o\`u ${\cal M}$ est un
fibr\'e holomorphe en droites (unique, \`a isomorphisme pr\`es) tel que
$c_1(K_X^{\vee}\otimes{\cal M}^{\otimes 2})=\pm c$. \hfill{\break}
3. ${\cal W}^g_X(c)$ est lisse. Soit $D$ le diviseur d'une section
non-triviale de ${\cal M}$. Alors ${\cal W}^g_X(c)$ a la dimension
pr\'esum\'ee si et seulement si
\hbox{$h^1({\cal O}_X(D)|_D)=0$}.
\end{thfr}

L'id\'ee de la d\'emonstration est d'identifier les classes d'isomorphie de
solutions de $(SW_0)$ avec les classes d'isomorphie de solutions d'une
\'equation de Vortex g\'en\'eralis\'ee [2], [7]. On utilise les m\'ethodes de
Bradlow pour montrer que celles-ci correspondent bijectivement aux sections
non-triviales de ${\cal M}$, modulo l'op\'eration naturelle du groupe $\C^*$.

Comme application, nous d\'emontrons le:
\begin{thfr} {\rm [8]}
Une surface complexe $X$ diff\'eomorphe \`a une surface rationnelle est
rationnelle.
\end{thfr}

La d\'emonstration consiste en trois parties:\\
i) Premi\`erement, soit $X_0$ une surface rationnelle. Il est bien connu [5]
que
$X_0$ admet une {\sl m\'etrique} de Hitchin , i.e.\ une m\'etrique
k\"ahlerienne $g_0$ telle que $c_1(K_{X_0})\cup[\omega_{g_0}]<0$. D'apr\`es le
th\'eor\`eme pr\'ec\'edent il s'ensuit que tous les espaces de modules
$W^{g_0}_X(c)$ sont vides, donc $n_c^{g_0}=0$ pour tout $c$ tel que $g_0$
est $c$-bonne.\\
ii) Supposons que ${\rm kod} X >0$ et que $X$ est l'\'eclatement en $k$
points de son mod\`ele minimal $X_{\min}$, et d\'esignons par
$\sigma:X\map X_{\min}$ la projection sur le mod\`ele minimal. On note
$E=\sum\limits_{i=1}^k E_i$ le diviseur exceptionnel, et pour
$I\subset\{1,\dots,k\}$ on pose $E_I:=\sum\limits_{i\in I} E_i$,
$\bar I:=\{1,\dots,k\}\setminus I$, $c_I:=2[E_I]-[K_X]$.

On fixe un diviseur ample $H_{\min}$ sur $X_{\min}$ et un nombre $n$
suffisamment grand. Soit $H_n:=\sigma^*(n H_{\min})-E$ le diviseur ample
correspondant sur $X$, et soit $g_n$ une m\'etrique k\"ahlerienne
telle que $[\omega_{g_n}]=[H_n]$. On d\'eduit facilement
$[\omega_{g_n}]\cdot c_I<0$ pour $n\gg 0$, et d'apr\`es le th\'eor\`eme
pr\'ec\'edent, $n_{c_I}^{g_n}\equiv 1$. \\
iii) Supposons qu'il existe un diff\'eomorphisme $f:X\map X_0$ compatible
avec les orientations, et consid\'erons l'image inverse $g:=f^*(g_0)$ d'une
m\'etrique de Hitchin sur $X$. On d\'eduit $n_{c_I}^g=0$ pour tout
$I\subset\{1,\dots,k\}$ tel que $g$ est $c$-bonne. On fait la remarque
suivante.\\
{\bf Remarque:} l'un des deux cas suivants a lieu: $g$ et $g_n$ sont
simultan\'ement \hbox{$c_I$-bonnes} et $n_{c_I}^g=n_{c_I}^{g_n}$, {\sl ou}
$g$ et $g_n$ sont simultan\'ement $c_{\bar I}$-bonnes et
\hbox{$n_{c_{\bar{I}}}^g=n_{c_{\bar I}}^{g_n}$}.  \\
En effet, soit $k:=\sigma^*([K_{\min}]_{\rm DR})$. On a $k^2\geq 0$ et
$k\ne 0$, donc une des deux composantes connexes de
${\cal K}$ est ${\cal K}_+:=\{u\in H^2_{\rm DR}(X)\ |\ u^2>0,u\cdot k>0\}$. On
v\'erifie facilement que $[\omega_{g_n}]\in{\cal K}_+$ et on choisit une forme
$g$-harmonique auto-duale $\omega_g$ telle que $[\omega_g]\in\KK_+$. Il suffit
maintenant de d\'emontrer que les demi-droites $\R_{>0}[\omega_{g}]$ et
$\R_{>0}[\omega_{g_n}]$ se trouvent du m\^eme c\^ot\'e de l'hyperplan $c_I$ ou
du m\^eme c\^ot\'e de l'hyperplan $c_{\bar I}$. Sinon, en \'ecrivant
$[\omega_{g}]=\sum\limits_{i=1}^k \lambda_i[E_i]+\sigma^*[\omega]$, avec
$[\omega]\in H^2_{\rm DR}(X_{\min})$ on d\'eduit:
$$-\sum\limits_{i\in I}\lambda_i+\sum\limits_{j\in\bar I}\lambda_j-
[\omega]\cdot k\geq 0 \ \ \ {\rm et}\ \ \
-\sum\limits_{j\in\bar I}\lambda_j+\sum\limits_{i\in I}\lambda_i-
[\omega]\cdot k\geq 0 .$$
Il en r\'esulte $[\omega]\cdot k\leq 0$, ce qui contradit le choix de
$[\omega]\in{\cal K}_+$, et ach\`eve la d\'emonstration du th\'eor\`eme.

Utilisant les r\'esultats obtenus dans [3], on conclut
\begin{cofr}
La dimension de Kodaira d'une surface complexe est un invariant de la
structure diff\'erentielle.
\end{cofr}
\parindent0cm
\vspace{0.5cm}
{{\bf R\'ef\'erences bibliographiques}}
\vspace{0.3cm}\small

[1] Barth, W., Peters, C., Van de Ven, A.: {\it Compact complex surfaces},
Springer Verlag, 1984.

[2] Bradlow, S. B.: {\it Vortices in holomorphic line bundles over closed
K\"ahler manifolds}, Comm.\ Math.\ Phys.\ 135,  1990, p.\ 1-17.

[3] Friedman, R., Morgan, J.W.: {\it Smooth 4-manifolds and Complex Surfaces},
Springer Verlag  3.\ Folge, Band 27, 1994.

[4]  Friedman, R., Qin, Z.: {\it On complex surfaces diffeomorphic to
rational surfaces}, Inventiones Math.\ (\`a para\^itre).

[5] Hitchin, N.: {\it  On the curvature of rational surfaces}, Proc.\ of
Symp.\ in Pure Math., Stanford, Vol.\ 27, 1975.

[6] Kronheimer, P., Mrowka, T.: {\it The genus of embedded surfaces in the
projective plane}, Preprint, Oxford, 1994.

[7] Okonek, Ch.; Teleman A.: {\it The Coupled Seiberg-Witten Equations,
Vortices, and Moduli Spaces of Stable Pairs}, Preprint, Z\"urich, Jan. 13-th,
1995.

[8] Okonek, Ch.; Teleman A.: {\it Seiberg-Witten Invariants and the Van De
Ven Conjecture}, Preprint, Z\"urich, February 8-th, 1995.

[9] Witten, E.: {\it Monopoles and four-manifolds}, Mathematical  Research
Letters 1, 1994, p.\ 769-796.

[10] Van de Ven, A,: {\it On the differentiable structure of certain
algebraic surfaces}, S\'em.\ Bourbaki ${\rm n}^o$ 667, Juin, 1986.
\vspace{0.1cm}\\
\\
Mathematisches Institut, Universit\"at Z\"urich,
Winterthurerst. 190,
CH-8057 Z\"urich.\hfill\break
e-mail:okonek@math.unizh.ch ; teleman@math.unizh.ch

\end{document}